\documentclass[12pt]{amsart}
\usepackage{amsmath,amssymb,latexsym,esint,cite,mathrsfs}
\usepackage{verbatim,wasysym}
\usepackage[left=2.6cm,right=2.6cm,top=3cm,bottom=3cm]{geometry}
\makeatletter \@addtoreset{equation}{section} \makeatother

\usepackage{graphicx}
\usepackage{subfigure}
\usepackage{graphicx,epic,eepic}

\usepackage[colorlinks=true,urlcolor=blue,
citecolor=red,linkcolor=blue,linktocpage,pdfpagelabels,
bookmarksnumbered,bookmarksopen]{hyperref}
\numberwithin{equation}{section}
\newtheorem{theorem}{Theorem}[section]
\newtheorem{lemma}[theorem]{Lemma}

\newtheorem{proposition}[theorem]{Proposition}
\newtheorem{corollary}[theorem]{Corollary}
\newtheorem{remark}[theorem]{Remark}
\numberwithin{equation}{section}

\allowdisplaybreaks

\begin{document}

\title[On a weighted fourth-order equation]
{Classification and non-degeneracy of positive radial solutions for a weighted fourth-order equation and its application}

\author[S. Deng]{Shengbing Deng}
\address{\noindent Shengbing Deng
\newline
School of Mathematics and Statistics, Southwest University,
Chongqing 400715, People's Republic of China}\email{shbdeng@swu.edu.cn}

\author[X. Tian]{Xingliang Tian$^{\ast}$}
\address{\noindent Xingliang Tian  \newline
School of Mathematics and Statistics, Southwest University,
Chongqing 400715, People's Republic of China.}\email{xltian@email.swu.edu.cn}

\thanks{$^{\ast}$ Corresponding author}



\allowdisplaybreaks
%\maketitle

\begin{abstract}
{\tiny This paper is devoted to radial solutions of the following weighted fourth-order equation
    \begin{equation*}
    \mathrm{div}(|x|^{\alpha}\nabla(\mathrm{div}(|x|^\alpha\nabla u)))=|u|^{2^{**}_{\alpha}-2}u\quad \mbox{in}\quad \mathbb{R}^N,
    \end{equation*}
where $N\geq 2$, $\frac{4-N}{2}<\alpha<2$ and $2^{**}_{\alpha}=\frac{2N}{N-4+2\alpha}$.  It is obvious that the solutions of above equation are invariant under the scaling $\lambda^{\frac{N-4+2\alpha}{2}}u(\lambda x)$  while they are not invariant under translation when $\alpha\neq 0$. We characterize all the solutions to the related linearized problem about radial solutions, and obtain the conclusion of that if $\alpha$ satisfies $(2-\alpha)(2N-2+\alpha)\neq4k(N-2+k)$ for all $k\in\mathbb{N}^+$ the radial solution is non-degenerate, otherwise there exist new solutions to the linearized problem that ``replace'' the ones due to the translations invariance.
As applications, firstly we investigate the remainder terms of some inequalities related to above equation. Then when $N\geq 5$ and $0<\alpha<2$, we establish a new type second-order Caffarelli-Kohn-Nirenberg inequality \begin{equation*}
    \int_{\mathbb{R}^N} |\mathrm{div}(|x|^\alpha\nabla u)|^2 \mathrm{d}x \geq C \left(\int_{\mathbb{R}^N}|u|^{2^{**}_{\alpha}} \mathrm{d}x\right)^{\frac{2}{2^{**}_{\alpha}}},\quad
    \mbox{for all}\quad u\in C^\infty_0(\mathbb{R}^N),
    \end{equation*}
    and in this case we consider a prescribed perturbation problem by using Lyapunov-Schmidt reduction.
    }
\end{abstract}

\vspace{3mm}

\maketitle

\section{{\bfseries Introduction}}\label{sectir}

    In this paper, we consider the following weighted fourth-order equation
    \begin{equation}\label{Pwb}
    \mathrm{div}(|x|^{\alpha}\nabla(\mathrm{div}(|x|^\alpha\nabla u)))=|u|^{2^{**}_{\alpha}-2}u\quad \mbox{in}\quad \mathbb{R}^N,
    \end{equation}
    where
    \begin{equation}\label{condpps}
    N\geq 2,\quad\frac{4-N}{2}<\alpha<2,\quad 2^{**}_{\alpha}=\frac{2N}{N-4+2\alpha}.
    \end{equation}
    This problem, for $\alpha\neq 0$, generalizes the well-known equation which involves the critical Sobolev exponent
    \begin{equation}\label{bce}
    \Delta^2 u=|u|^{2^{**}-2}u, \quad \mbox{in}\quad \mathbb{R}^N,
    \end{equation}
    where $2^{**}=\frac{2N}{N-4}$. Smooth positive solutions to (\ref{bce}) have been completely classified by C.-S. Lin \cite{Li98}, where the authors proved that they are given by
    \begin{equation*}
    U[z,\lambda](x)=[(N-4)(N-2)N(N+2)]^{\frac{N-4}{8}}\left(\frac{\lambda}{1+\lambda^2|x-z|^2}\right)^{\frac{N-4}{2}},\quad \lambda>0,\quad z\in\mathbb{R}^N,
    \end{equation*}
    and they are extremal functions for the well-known Sobolev inequality
    \begin{equation}\label{bcesi}
    \int_{\mathbb{R}^N}|\Delta u|^2 \mathrm{d}x\geq \mathcal{S}_0\left(\int_{\mathbb{R}^N}|u|^{2^{**}} \mathrm{d}x\right)^\frac{2}{2^{**}},\quad \mbox{for all}\quad u\in \mathcal{D}^{2,2}_0(\mathbb{R}^N),
    \end{equation}
    for some constant $\mathcal{S}_0$, where $\mathcal{D}^{2,2}_0(\mathbb{R}^N):=\{u\in L^{2^{**}}(\mathbb{R}^N): \Delta u\in L^2(\mathbb{R}^N)\}$.

    The presence of the term $|x|^\alpha$ in equation (\ref{Pwb}) drastically changes the problem. For this kind of operators it is not possible to apply the moving plane method anymore (to get the radial symmetry around some point), and indeed non-radial solutions may appear.
    Problem (\ref{Pwb}) is motivated by Gladiali, Grossi and Neves \cite{GGN13}, where the authors considered the second-order H\'{e}non problem
    \begin{equation}\label{Phs}
    -\Delta u=(N+l)(N-2)|x|^l u^{\frac{N+2+2l}{N-2}},\quad u>0 \quad \mbox{in}\quad \mathbb{R}^N,
    \end{equation}
    where $N\geq 3$ and $l>0$. This problem generalizes the well-known critical Sobolev inequality.
    Firstly, the authors gave the classification of radial solutions $V^{(\lambda)}_{l}$ in classical Sobolev space $\mathcal{D}^{1,2}_0(\mathbb{R}^N):=\{u\in L^{\frac{2N}{N-2}}(\mathbb{R}^N): \nabla u\in L^2(\mathbb{R}^N)\}$ of equation (\ref{Phs}), where $V^{(\lambda)}_{l}(x)=\lambda^\frac{N-2}{2}V_l(\lambda x)$ and $V_l(x)=(1+|x|^{2+l})^{-(N-2)/(2+l)}$. Furthermore, they characterized all the solutions to the linearized problem related to (\ref{Phs}) at function $V_l$, that is,
    \begin{equation}\label{Pwhls}
    -\Delta v=(N+l)(N+2+2l)|x|^l V_l^{\frac{4+2l}{N-2}}v \quad \mbox{in}\quad \mathbb{R}^N, \quad v\in \mathcal{D}^{1,2}_0(\mathbb{R}^N).
    \end{equation}

\vskip0.25cm

    \noindent{\bf Theorem~A.} \cite[Theorem 1.3]{GGN13} {\it Let $l\geq 0$.  If $l>0$ is not an even integer, then the space of solutions of (\ref{Pwhls}) has dimension one and is spanned by
    \begin{equation*}
    X_0(x)=\frac{1-|x|^{2+l}}{(1+|x|^{2+l})^\frac{N+l}{2+l}},
    \end{equation*}
    where $X_0\sim\frac{\partial V^{(\lambda)}_{l}}{\partial \lambda}|_{\lambda=1}$. If $l=2(k-1)$ for some $k\in\mathbb{N}^+$, then the space of solutions of (\ref{Pwhls}) has dimension $1+\frac{(N+2k-2)(N+k-3)!}{(N-2)!k!}$ and is spanned by
    \begin{equation*}
    X_0(x)=\frac{1-|x|^{2+l}}{(1+|x|^{2+l})^\frac{N+l}{2+l}},\quad X_{k,i}(x)=\frac{|x|^k\Psi_{k,i}(x)}{(1+|x|^{2+l})^\frac{N+l}{2+l}},
    \end{equation*}
    where $\{\Psi_{k,i}\}$, $i=1,\ldots,\frac{(N+2k-2)(N+k-3)!}{(N-2)!k!}$, form a basis of $\mathbb{Y}_k(\mathbb{R}^N)$, the space of all homogeneous harmonic polynomials of degree $k$ in $\mathbb{R}^N$.
    }

\vskip0.25cm

\noindent As stated in \cite[p.4]{GGN13}, Theorem A highlights the new phenomenon that if $l$ is an even integer then there exist new solutions to (\ref{Pwhls}) that ``replace'' the ones due to the translations invariance. Furthermore, the authors in \cite{GGN13} found a bifurcating phenomenon at $(l,V_l)$ when $l$ is an even integer that is nonradial solutions appear.

    Then let us mention the following weight fourth-order elliptic equation
    \begin{equation}\label{Pl}
    \Delta(|x|^{\alpha}\Delta u)=|x|^l u^p,\quad u\geq 0 \quad \mbox{in}\quad \mathbb{R}^N,
    \end{equation}
    where $N\geq 5$, $p>1$ and $4-N<\alpha<\min\{N,l+4\}$.
    Define
    \begin{equation*}%\label{defps}
    p_s:=\frac{N+4+2l-\alpha}{N-4+\alpha}.
    \end{equation*}
    Guo et al. \cite{GWW20}  obtained Liouville type result, that is, if $u\in C^4(\mathbb{R}^N\backslash\{0\}) \cap C^0(\mathbb{R}^N)$ with $|x|^{\alpha}\Delta u\in C^0(\mathbb{R}^N)$ is a nonnegative radial solution to (\ref{Pl}), then $u\equiv 0$ in $\mathbb{R}^N$ provided $1<p<p_s$. Huang and Wang in \cite{HW20} got the partial classifications of positive radial solutions of equation (\ref{Pl}) with $p=p_s$ for some special cases, see also \cite{Ya21} for more general case.
    We have established the analogous conclusion as \cite[Theorem 1.3]{GGN13} related to the linearized problem for equation (\ref{Pl}) with $p=p_s$ for $l=-\alpha$ or $l=\frac{4(\alpha-2)(N-2)}{N-\alpha}-\alpha$ respectively, see \cite{DGT23-jde,DT22} for details. Recently, Guan and Guo \cite{GG22} considered the existence and regularity of positive solutions for a new type weighted fourth-order problem.

    Therefore, it is natural to consider that could we establish the analogous conclusion as
    \cite[Theorem 1.3]{GGN13} related to some other kinds of fourth-order problems, that is, in this paper we consider problem (\ref{Pwb}). Define
    \begin{equation}\label{defd2a}
    \mathcal{D}^{2,2}_\alpha(\mathbb{R}^N):=cl\left\{u\in C^\infty_c(\mathbb{R}^N): \int_{\mathbb{R}^N} |\mathrm{div}(|x|^\alpha\nabla u)|^2 \mathrm{d}x<\infty\right\},
    \end{equation}
    with respect to the inner product
    \[
    \langle u,v\rangle_\alpha
    := \int_{\mathbb{R}^N} \mathrm{div}(|x|^\alpha\nabla u)\mathrm{div}(|x|^\alpha\nabla v) \mathrm{d}x
    \]
    and the norm $\|u\|_{\mathcal{D}^{2,2}_\alpha(\mathbb{R}^N)}:=\langle u,u\rangle_\alpha^{1/2}$. Firstly, by using the Emden-Folwer translation as in \cite{HW20}, we obtain the uniqueness of radial solutions of \eqref{Pwb}.

    \begin{theorem}\label{thmpwh}
    Suppose $N\geq 2$ and $\frac{4-N}{2}<\alpha<2$. Then problem \eqref{Pwb} in $\mathcal{D}^{2,2}_\alpha(\mathbb{R}^N)$ has a unique (up to scalings and change of sign)
    radial nontrivial solution of the form $\pm U_{\lambda,\alpha}(x)$ for some $\lambda>0$, where
    \begin{equation}\label{defula}
    U_{\lambda,\alpha}(x)=\frac{C_{N,\alpha}\lambda^{\frac{N-4+2\alpha}{2}}}{(1+\lambda^{2-\alpha}|x|^{2-\alpha})^{\frac{N-4+2\alpha}{2-\alpha}}},
    \end{equation}
    with $C_{N,\alpha}=\left[(N-4+2\alpha)(N-2+\alpha)N(N+2-\alpha)\right]^{\frac{N-4+2\alpha}{8-4\alpha}}$.
    \end{theorem}

    Then as a corollary of Theorem \ref{thmpwh}, we directly obtain the following result, which in the radial case, extends the classical inequality (\ref{bcesi}).

    \begin{corollary}\label{thmPbcb}
    Suppose $N\geq 2$ and $\frac{4-N}{2}<\alpha<2$. Let $u\in \mathcal{D}^{2,2}_\alpha(\mathbb{R}^N)$ be a radial function. Then we have that,
    \begin{equation}\label{Ppbcm}
    \int_{\mathbb{R}^N} |\mathrm{div}(|x|^\alpha\nabla u)|^2 \mathrm{d}x \geq S^{rad}_{N,\alpha}\left(\int_{\mathbb{R}^N}|u|^{2^{**}_{\alpha}} \mathrm{d}x \right)^{\frac{2}{2^{**}_{\alpha}}},
    \end{equation}
    for some positive constant $S^{rad}_{N,\alpha}$. The explicit form of the best constant in (\ref{Ppbcm}) is
    \begin{equation*}
    S^{rad}_{N,\alpha}
    =\left(\frac{2-\alpha}{2}\right)^{4-\frac{4-2\alpha}{N}}
    \left(\frac{2\pi^{\frac{N}{2}}}{\Gamma(\frac{N}{2})}
    \right)^{\frac{4-2\alpha}{N}}C\left(\frac{2N}{2-\alpha}\right),
    \end{equation*}
    where $C(M)=(M-4)(M-2)M(M+2)\left[\Gamma^2(\frac{M}{2})/(2\Gamma(M))\right]^{\frac{4}{M}}$.
    Moreover the extremal radial functions which achieve $S^{rad}_{N,\alpha}$ (equality holds in \eqref{Ppbcm}) are given by
    \begin{equation*}%\label{pbcm}
    V_{\lambda,\alpha}(x)=\frac{A\lambda^{\frac{N-4+2\alpha}{2}}}
    {(1+\lambda^{2-\alpha}|x|^{2-\alpha}
    )^{\frac{N-4+2\alpha}{2-\alpha}}},\quad \mbox{for some}\ A\in\mathbb{R},\ \lambda>0.
    \end{equation*}
    \end{corollary}

    Inspired by \cite{GGN13}, then we concern the linearized problem related to \eqref{Pwb} at the function $U_{1,\alpha}$. This leads to study the problem
    \begin{equation}\label{Pwhl}
    \mathrm{div}(|x|^{\alpha}\nabla(\mathrm{div}(|x|^\alpha\nabla v)))=(2^{**}_{\alpha}-1)U_{1,\alpha}^{2^{**}_{\alpha}-2}v \quad \mbox{in}\quad \mathbb{R}^N, \quad v\in \mathcal{D}^{2,2}_\alpha(\mathbb{R}^N).
    \end{equation}
    Next theorem characterizes all the solutions to (\ref{Pwhl}).

    \begin{theorem}\label{thmpwhl}
    Suppose $N\geq 2$ and $\frac{4-N}{2}<\alpha<2$. If
    \begin{equation}\label{ncqft}
    \left(\frac{2-\alpha}{2}\right)^2\left(\frac{2N}{2-\alpha}-1\right)=k(N-2+k),\quad \mbox{for some}\quad k\in\mathbb{N}^+,
    \end{equation}
    then the space of solutions of (\ref{Pwhl}) has dimension $1+\frac{(N+2k-2)(N+k-3)!}{(N-2)!k!}$ and is spanned by
    \begin{equation}\label{defaezki}
    Z_{0}(x)=\frac{1-|x|^{2-\alpha}}{(1+|x|^{2-\alpha})^\frac{N-2+\alpha}{2-\alpha}},\quad Z_{k,i}(x)=\frac{|x|^{\frac{2-\alpha}{2}}\Psi_{k,i}(x)}{(1+|x|^{2-\alpha})^\frac{N-2+\alpha}{2-\alpha}},
    \end{equation}
    where $\{\Psi_{k,i}\}$, $i=1,\ldots,\frac{(N+2k-2)(N+k-3)!}{(N-2)!k!}$, form a basis of $\mathbb{Y}_k(\mathbb{R}^N)$, the space of all homogeneous harmonic polynomials of degree $k$ in $\mathbb{R}^N$. Otherwise, that is, if
     \begin{equation}\label{ncqftb}
    (2-\alpha)(2N-2+\alpha)\neq4k(N-2+k),\quad \mbox{for all}\quad k\in\mathbb{N}^+,
    \end{equation}
    the solution $U_{1,\alpha}$ of equation (\ref{Pwb}) is non-degenerate in the sense that the space of solutions of (\ref{Pwhl}) has dimension one and is spanned by $Z_0$.
    \end{theorem}

    \begin{remark}\label{rem:exf}
    The key step of the proof for Theorem \ref{thmpwh} is making use of Endem-Folwer translation (see (\ref{l-et})), while the key step of the proof for Theorem \ref{thmpwhl} is the change of variable $r\mapsto r^{\frac{2-\alpha}{2}}$ which changes the dimension $N$ into $M=\frac{2N}{2-\alpha}>4$.

    We note that in the case $\alpha=0$ we get $k=1$ and one gets back the well known result for the equation involving the critical Sobolev exponent, see \cite{Fe02,LW00}. We observe that for $\alpha\neq 0$ the solutions of problem (\ref{Pwb}) are invariant for dilations but not for translations. Theorem \ref{thmpwhl} highlights the new phenomenon that if $\alpha$ satisfies (\ref{ncqft}) then there exist new solutions to (\ref{Pwhl}) that ``replace'' the ones due to the translations invariance as in Theorem A.

    \end{remark}

    As applications of Corollary \ref{thmPbcb} and Theorem \ref{thmpwhl}, inspired by \cite{BrE85}, firstly we consider the remainder terms of inequality (\ref{Ppbcm}) in radial space
    \begin{equation}\label{defd2ar}
    \mathcal{D}^{2,2}_{\alpha,rad}(\mathbb{R}^N):=\{u\in \mathcal{D}^{2,2}_\alpha(\mathbb{R}^N): u(x)=u(|x|)\},
    \end{equation}
    where $\mathcal{D}^{2,2}_\alpha(\mathbb{R}^N)$ is defined as in (\ref{defd2a}), as an analogous result to \cite{LW00}.

    \begin{theorem}\label{thmprt}
    Suppose $N\geq 2$ and $\frac{4-N}{2}<\alpha<2$. Then there exists constant $B=B(N,\alpha)>0$ such that for every $u\in \mathcal{D}^{2,2}_{\alpha,rad}(\mathbb{R}^N)$, it holds that
    \[
    \int_{\mathbb{R}^N} |\mathrm{div}(|x|^\alpha\nabla u)|^2 \mathrm{d}x
    -S^{rad}_{N,\alpha}\left(\int_{\mathbb{R}^N}|u|^{2^{**}_{\alpha}} \mathrm{d}x \right)^{\frac{2}{2^{**}_{\alpha}}}
    \geq B\mathrm{dist}(u,\mathcal{M}_\alpha)^2,
    \]
    where $\mathcal{M}_\alpha=\{cU_{\lambda,\alpha}: c\in\mathbb{R}, \lambda>0\}$ is set of extremal functions for \eqref{Ppbcm}, and $\mathrm{dist}(u,\mathcal{M}_\alpha):=\inf\limits_{c\in\mathbb{R}, \lambda>0}\|u-cU_{\lambda,\alpha}\|_{\mathcal{D}^{2,2}_\alpha(\mathbb{R}^N)}$.
    \end{theorem}

    The second thing we want to study is to construct solutions by using the Lyapunov-Schmidt argument, enlightened by \cite{AGP99,Fe02}
    (and also \cite[Sections 3 and 4]{FS03}). When $\alpha=0$ and $N\geq 5$, Felli \cite{Fe02} considered the following scalar curvature problem
    \begin{equation}\label{Pwhpu}
    \Delta^2u=(1+\varepsilon K(x))|u|^{\frac{8}{N-4}}u\quad \mbox{in}\quad \mathbb{R}^N,
    \end{equation}
    under some suitable assumptions on $K$, the author proved \eqref{Pwhpu} admits a solution for any $\varepsilon>0$ small. In \cite{DT22}, we have considered a similar weighted fourth-order perturbation  problem in radial space when $\alpha$ is not an even integer. In present paper, we will remove the radial assumption when $0<\alpha<2$ and $N\geq 5$.

    Note that $2^{**}_{\alpha}\in (2,\frac{2N}{N-4})$ when $0<\alpha<2$ and $N\geq 5$. In this case, let us take the same arguments as those in \cite[Section 2]{GG22} (which deals with it in bounded domain, but it also holds in the whole space), then we deduce that for all $u\in C^\infty_0(\mathbb{R}^N)$ the functional
    \[
    \int_{\mathbb{R}^N} |\mathrm{div}(|x|^\alpha\nabla u)|^2 \mathrm{d}x
    \]
    is equivalent to
    \[
    \int_{\mathbb{R}^N} |x|^{2\alpha}|\Delta u|^2 \mathrm{d}x.
    \]
    For the convenience of readers, we give the proof in Proposition \ref{propneq}.
    Therefore, by the embedding theorem as in \cite[Lemma A.1]{MS14}, we deduce that $\mathcal{D}^{2,2}_{\alpha}(\mathbb{R}^N)\hookrightarrow
    L^{2^{**}_{\alpha}}(\mathbb{R}^N)$ continuously, that is, there exists a constant $S_{N,\alpha}>0$ such that
    \begin{equation}\label{ckncc}
    \int_{\mathbb{R}^N} |\mathrm{div}(|x|^\alpha\nabla u)|^2 \mathrm{d}x \geq S_{N,\alpha} \left(\int_{\mathbb{R}^N}|u|^{2^{**}_{\alpha}} \mathrm{d}x\right)^{\frac{2}{2^{**}_{\alpha}}},\quad
    \mbox{for all}\quad u\in \mathcal{D}^{2,2}_{\alpha}(\mathbb{R}^N).
    \end{equation}
    This is a new type second-order Caffarelli-Kohn-Nirenberg inequality, comparing with \cite{Li86}. Whether the best constant $S_{N,\alpha}$ equals to $S^{rad}_{N,\alpha}$ is not known yet, because the symmetric decreasing rearrangement method does not make sense in the case $\alpha>0$. This is a challenge problem.

    Now, by using above inequality and Theorem \ref{thmpwhl}, we will establish sufficient conditions on a prescribed weighted $h(x)$ on $\mathbb{R}^N$ which guarantee the existence of solutions to the following perturbation problem
    \begin{equation}\label{Pwhp}
    \mathrm{div}(|x|^{\alpha}\nabla(\mathrm{div}(|x|^\alpha\nabla u)))=(1+\varepsilon h(x))u^{2^{**}_{\alpha}-1},\quad u>0 \quad \mbox{in}\quad \mathbb{R}^N, \quad u\in \mathcal{D}^{2,2}_{\alpha}(\mathbb{R}^N).
    \end{equation}

    \begin{theorem}\label{thmpwhp}
    Assume that $N\geq 5$, $0<\alpha<2$, and $h\in L^\infty(\mathbb{R}^N)\cap C(\mathbb{R}^N)$. If $\lim\limits_{|x|\to 0}h(x)=\lim\limits_{|x|\to \infty}h(x)=0$, then for $|\varepsilon|$ small enough, (\ref{Pwhp}) has at least one solution of the form
    \[
    u_\varepsilon(x)=U_{\lambda,\alpha}(x)+\omega(\lambda,\varepsilon),
    \]
    where $\lambda=\lambda(\varepsilon)$ and $\|\omega(\lambda,\varepsilon)\|_{\mathcal{D}^{2,2}_\alpha(\mathbb{R}^N)}
    =O(\varepsilon)$.
    \end{theorem}

    \begin{remark}\label{rem:pp}
    If $\alpha$ is an even integer, we can refer to \cite{Al18,EGPV21,GG15} which deal with the second-order Hardy or H\'{e}non problems in bounded domain.
    \end{remark}

    The paper is organized as follows: In Section \ref{sectpmr} we deduce the uniqueness of radial solutions of \eqref{Pwb} and give the optimizers of inequality (\ref{Ppbcm}). Section \ref{sectlp} is devoted to characterizing all solutions to the linearized equation (\ref{Pwhl}). In Section \ref{sect:rt}, we study the stability of inequality (\ref{Ppbcm}) and give the proof of Theorem \ref{thmprt}. In Section \ref{sectprp} we investigate the existence of solutions to the perturbed equation (\ref{Pwhp}) by using finite-dimensional reduction method and prove Theorem \ref{thmpwhp}.

\section{{\bfseries Optimizers of inequality (\ref{Ppbcm})}}\label{sectpmr}

Firstly, we give the classification of radial solutions of \eqref{Pwb}.
In order to obtain this, we will use the classical Emden-Fowler transformation.
    Following the work of Huang and Wang \cite{HW20}, let $\varphi, \psi\in C^2(\mathbb{R})$ be defined by
    \begin{align}\label{l-et}
    u(|x|)=|x|^{-\kappa_1}\varphi(t),\quad v(|x|)=|x|^{-\kappa_2}\phi(t),\quad \mbox{with}\ t=-\ln |x|,
    \end{align}
    where $-\mathrm{div} (|x|^\alpha \nabla u)=v$, and
    \begin{align*}
    \kappa_1=\frac{N}{2^{**}_{\alpha}}
    =\frac{N-4+2\alpha}{2},\quad
    \kappa_2=2-\alpha+\kappa_1=\frac{N}{2}.
    \end{align*}
    Since $\alpha>(4-N)/2$ then $\kappa_1+\kappa_2=N-2+\alpha>0$.

    Denote the constants
    \begin{align*}
    \mathcal{A}:& = \frac{N-2+\alpha}{2}-\kappa_1
    =-\frac{N-2+\alpha}{2}+\kappa_2
    =\frac{2-\alpha}{2}, \\
    \mathcal{B}:& = \kappa_1(N-2+\alpha-\kappa_1)=\kappa_1\kappa_2
    =\frac{N(N-4+2\alpha)}{4}.
    \end{align*}

\medskip

\noindent {\bf \em Proof of Theorem \ref{thmpwh}.}
Direct computation shows that a radial function $u\in \mathcal{D}^{2,2}_{\alpha}(\mathbb{R}^N)$ solves \eqref{Pwb} if and only if the functions $\varphi, \phi$ satisfy
    \begin{eqnarray*}%\label{Pwh}
    \left\{ \arraycolsep=1.5pt
       \begin{array}{ll}
        -\varphi''+2\mathcal{A}\varphi'+\mathcal{B}\varphi=\phi\quad \mbox{in}\ \mathbb{R},\\[2mm]
        -\phi''-2\mathcal{A}\phi'+\mathcal{B}\phi
        =|\varphi|^{2^{**}_{\alpha}-2}\varphi\quad \mbox{in}\ \mathbb{R},
        \end{array}
    \right.
    \end{eqnarray*}
    that is, $\varphi(t)$ satisfies the following fourth-order ordinary differential equation
    \begin{equation}\label{Pwht}
    \varphi^{(4)}-K_2\varphi''+K_0\varphi
    =|\varphi|^{2^{**}_{\alpha}-2}\varphi\quad \mbox{in}\  \mathbb{R},\quad \varphi\in H^2(\mathbb{R}),
    \end{equation}
    with the constants
    \begin{align*}
    K_2= \kappa_1^2+\kappa_2^2
    =\frac{(N+\alpha-2)^2+(2-\alpha)^2}{2}, \quad
    K_0= \kappa_1^2\kappa_2^2
    =\frac{N^2(N-4+2\alpha)^2}{16}.
    \end{align*}
    Here $H^2(\mathbb{R})$ denotes the completion of $C^\infty_0(\mathbb{R})$ with respect to the norm
    \[
    \|\varphi\|_{H^2(\mathbb{R})}
    =\left(\int_{\mathbb{R}}(|\varphi_{tt}|^2+K_2|\varphi_t|^2
    +K_0|\varphi|^2)\mathrm{d}t\right)^{1/2}.
    \]
    Clearly, $K_2^2\geq 4 K_0$. Then since $2^{**}_{\alpha}>2$, applying
    \cite[Theorem 2.2 ]{BM12} directly, we could get the existence and uniqueness (up to translations, inversion $t\mapsto -t$ and change of sign) of smooth solution to equation \eqref{Pwht}.

    Indeed, we can set
    \begin{equation}\label{vtd}
    \varphi(t)=\mathcal{C}(\cosh \nu t)^m,
    \end{equation}
    where $\mathcal{C}, \nu, m$ are to be determined later. Here $\cosh s=(e^s+e^{-s})/2$. As in \cite{HW20}, it is easy to verify that
    \begin{align}\label{mnc}
    m= & -\frac{4}{2^{**}_{\alpha}-2}=\frac{N}{2-\alpha},\quad \nu=\left(\frac{K_2}{m^2+(m-2)^2}\right)^{\frac{1}{2}}
    =\frac{\alpha-2}{2},\nonumber\\
    \mathcal{C}= & \left[m(m-1)(m-2)(m-3)\nu^4\right]
    ^{\frac{1}{2^{**}_{\alpha}-2}}.
    \end{align}
    Now, let us return to the original problem. The change of $u(|x|)=|x|^{-\frac{N-4+2\alpha}{2}}\varphi(t)$ with $t=-\ln |x|$, then $\varphi(t)=\mathcal{C}(\cosh \nu t)^m$ indicates
    \[
    u(|x|)=\frac{\mathcal{C}}{2^m}|x|^{-(\frac{N-4+2\alpha}{2}+m\nu)}
    (1+|x|^{2\nu})^m,
    \]
    where $m, \nu, \mathcal{C}$ are given as in \eqref{mnc}. Thus we obtain
    \begin{equation*}%\label{defulae}
    u(x)=\frac{C_{N,\alpha}}{(1+|x|^{2-\alpha})^{\frac{N-4+2\alpha}{2-\alpha}}},
    \end{equation*}
    where $C_{N,\alpha}=\left[(N-4+\alpha)(N-2)(N-\alpha)(N+2-2\alpha)
    \right]^{\frac{N-4+2\alpha}{8-4\alpha}}$. Note that the invariance of $\varphi$ up to translations and inversion $t\mapsto -t$ indicates invariance of $u$ up to scalings, thus \eqref{Pwb} admits a unique (up to scalings and change of sign)
    radial solution of the form $\pm U_{\lambda,\alpha}(x)$ for some $\lambda>0$, where
    \begin{equation*}
    U_{\lambda,\alpha}(x)=\frac{C_{N,\alpha}\lambda^{\frac{N-4+2\alpha}{2}}}
    {(1+\lambda^{2-\alpha}|x|^{2-\alpha})^{\frac{N-4+2\alpha}{2-\alpha}}}.
    \end{equation*}
    Now, the proof of Theorem \ref{thmpwh} is completed.
    \qed

\medskip

    Corollary \ref{thmPbcb} can be done directly as those in \cite{HW20} by using the Emden-Fowler transformation shown as in \eqref{l-et}. %Here we will use different method that transforms the dimension $N$ into $M=\frac{2N}{2-\alpha}>2$ which was first introduced in \cite{GGN13} (see also \cite{DT22}), that is, we make the change of variable $r\mapsto r^{\frac{2-\alpha}{2}}$, related to Sobolev inequality to investigate the sharp constants and optimizers of inequality (\ref{Ppbcm}).
    Recall $\mathcal{\mathcal{D}}^{2,2}_{\alpha,rad}(\mathbb{R}^N)$ defined as in (\ref{defd2ar}), then the best constant in (\ref{Ppbcm}) can be defined by
    \begin{equation}\label{defbcs}
    S^{rad}_{N,\alpha}:=\inf_{u\in \mathcal{\mathcal{D}}^{2,2}_{\alpha,rad}(\mathbb{R}^N)\backslash\{0\}}
    \frac{\int_{\mathbb{R}^N} |{\rm div}(|x|^\alpha\nabla u)|^2 \mathrm{d}x}{\left(\int_{\mathbb{R}^N}|u|^{2^{**}_{\alpha}} \mathrm{d}x \right)^{\frac{2}{2^{**}_{\alpha}}}}.
    \end{equation}

\medskip

\noindent {\bf \em Proof of Corollary \ref{thmPbcb}.}
We follow the arguments as those in the proof of \cite[Theorem A.1]{GGN13}. Let $u\in \mathcal{\mathcal{D}}^{2,2}_{\alpha,rad}(\mathbb{R}^N)$. Making the change that $v(s)=u(r)$ with $|x|=r=s^{1/q}$ where $q>0$ will be given later, then we have
    \begin{equation*}
    \begin{split}
     \int_{\mathbb{R}^N} |{\rm div}(|x|^\alpha\nabla u)|^2 \mathrm{d}x
    & = \omega_{N-1}\int^\infty_0 \left[u''(r)+\frac{N+\alpha-1}{r}u'(r)\right]^2r^{N+2\alpha-1}
    \mathrm{d}r \\
    & = \omega_{N-1}q^{3}\int^\infty_0\left[v''(s)
    +\frac{(N-2+\alpha+q)/q}{s}v'(s)\right]^2
    s^{(N-4+2\alpha+3q)/q}\mathrm{d}s,
    \end{split}
    \end{equation*}
    where $\omega_{N-1}=2\pi^{\frac{N}{2}}/\Gamma(\frac{N}{2})$  is the surface area for unit ball of $\mathbb{R}^N$.
    In order to make use of the Sobolev inequality of radial case established by de Oliveira and Silva \cite{dS23}, we need $(N-2+\alpha+q)/q=(N-4+2\alpha+3q)/q$ which requires $q=\frac{2-\alpha}{2}.$ Now, we set
    \begin{equation}\label{defm}
    M:=\frac{N-2+\alpha+q}{q+1}=\frac{2N}{2-\alpha}>4,
    \end{equation}
    which implies
    \[
    \frac{2M}{M-4}=\frac{2N}{N-4+2\alpha}=2^{**}_{\alpha},\quad qM=N,
    \]
    therefore
    \begin{equation*}
    \begin{split}
    \int^\infty_0 \left[u''(r)+\frac{N+\alpha-1}{r}u'(r)\right]^2r^{N+2\alpha-1}
    \mathrm{d}r
    = q^{3}\int^\infty_0\left[v''(s)+\frac{M-1}{s}v'(s)\right]^2
    s^{M-1}\mathrm{d}s,
    \end{split}
    \end{equation*}
    and
    \begin{equation*}
    \begin{split}
    \int^\infty_0 |u(r)|^{2^{**}_{\alpha}}r^{N-1}
    \mathrm{d}r
    =\int^\infty_0
    |u(r)|^{\frac{2M}{M-4}}r^{qM-1}
    \mathrm{d}r
    =q^{-1}\int^\infty_0|v(s)|^{\frac{2M}{M-4}}s^{M-1}
    \mathrm{d}s.
    \end{split}
    \end{equation*}
    To sum up,
    \begin{equation}\label{defbcst}
    S^{rad}_{N,\alpha}=q^{4-\frac{4}{M}}\omega^{\frac{4}{M}}_{N-1}C(M),
    \end{equation}
    where
    \begin{equation}\label{defbcscm}
    C(M)=\inf_{v\in \mathcal{\mathcal{D}}^{2,2}_{\infty}(M)\backslash\{0\}}
    \frac{\int^\infty_0\left[v''(s)+\frac{M-1}{s}v'(s)\right]^2
    s^{M-1}\mathrm{d}s}
    {\left(\int^\infty_0|v(s)|^{\frac{2M}{M-4}}s^{M-1}
    \mathrm{d}s \right)^{\frac{M-4}{M}}}.
    \end{equation}
    Here $\mathcal{\mathcal{D}}^{2,2}_{\infty}(M)$ denotes the completion of $C^\infty_0((0,\infty))$ with respect to the norm
    \[
    \|v\|^2_{\mathcal{\mathcal{D}}^{2,2}_{\infty}(M)}=\int^\infty_0\left[v''(s)+\frac{M-1}{s}v'(s)\right]^2
    s^{M-1}\mathrm{d}s,
    \]
    see \cite{dS23}. In fact, the authors in \cite{dS23} proved that $C(M)>0$ which can be achieved, and the minimizers are solutions (up to some suitable multiplications) of
    \begin{equation}\label{les}
    \Delta_s^2v=|v|^{\frac{8}{M-4}}v\quad \mbox{in}\ (0,\infty),\quad v\in \mathcal{\mathcal{D}}^{2,2}_{\infty}(M),
    \end{equation}
    where $\Delta_s=\frac{\partial^2}{\partial s^2}+\frac{M-1}{s}\frac{\partial}{\partial s}$. Theorem \ref{thmpwh} indicates that \eqref{les} admits a unique (up to scalings and change of sign) solution of the form $v(s)=\pm U_{\lambda,\alpha}(s)$ given as in \eqref{defula} which implies $cU_{\lambda,\alpha}$ is the unique minimizer for $C(M)$ for some $c\in\mathbb{R}\setminus\{0\}$ and $\lambda>0$. Therefore, putting $U_{1,\alpha}$ into \eqref{defbcscm} as a test function, we can directly obtain
    \begin{align*}
    C(M) =(M-4)(M-2)M(M+2)
    \left(\frac{\Gamma^2\left(\frac{M}{2}\right)}
    {2\Gamma(M)}\right)^{\frac{4}{M}},
    \end{align*}
    where $\Gamma$ denotes the classical Gamma function.
    Then turning back to \eqref{defbcst}, we have
    \begin{equation*}
    S^{rad}_{N,\alpha}
    = \left(\frac{2-\alpha}{2}\right)^{4-\frac{4-2\alpha}{N}}
    \left(\frac{2\pi^{\frac{N}{2}}}{\Gamma(\frac{N}{2})}
    \right)^{\frac{4-2\alpha}{N}}C\left(\frac{2N}{2-\alpha}\right),
    \end{equation*}
    and it is achieved if and only if by
    \begin{equation*}%\label{defvla}
    V_{\lambda,\alpha}(x)=\frac{A\lambda^{\frac{N-4+2\alpha}{2}}}{(1+\lambda^{2-\alpha}|x|^{2-\alpha})^{\frac{N-4+2\alpha}{2-\alpha}}},
    \end{equation*}
    for some $A\in\mathbb{R}\backslash\{0\}$ and $\lambda>0$.
    The proof of Corollary \ref{thmPbcb} is now completed.
    \qed

\medskip

\section{{\bfseries Classification of linearized problem}}\label{sectlp}

By using the standard spherical decomposition and taking the change of variable $v(s)=u(r)$ with $s=r^{\frac{2-\alpha}{2}}$, we can characterize all solutions to the linearized problem (\ref{Pwhl}).

\medskip

\noindent {\bf \em Proof of Theorem \ref{thmpwhl}.}
    We follow the arguments as those in the proof of \cite[Theorem 1.3]{GGN13}.
    We decompose the fourth-order equation (\ref{Pwhl}) into a system of two second-order equations. Set
    \begin{equation}\label{PpwhlW}
    \begin{split}
    -\mathrm{div}(|x|^{\alpha}\nabla v)=W,
    \end{split}
    \end{equation}
    then problem (\ref{Pwhl}) is equivalent to the following system:
    \begin{eqnarray}\label{Pwhlp}
    \left\{ \arraycolsep=1.5pt
       \begin{array}{ll}
        -|x|^{\alpha}\Delta v-\alpha|x|^{\alpha-2}(x\cdot\nabla v)=W \quad \mbox{in}\  \mathbb{R}^N,\\[3mm]
        -|x|^{\alpha}\Delta W-\alpha|x|^{\alpha-2}(x\cdot\nabla W)
        =\frac{(2^{**}_{\alpha}-1)C_{N,\alpha}^{2^{**}_{\alpha}-2}}
        {(1+|x|^{2-\alpha})^4}v\quad \mbox{in}\  \mathbb{R}^N,
        \end{array}
    \right.
    \end{eqnarray}
    in $v\in\mathcal{D}^{2,2}_\alpha(\mathbb{R}^N)$.

    Firstly, we decompose $v$ and $W$ as follows:
    \begin{equation}\label{defvd}
    v(r,\theta)=\sum^{\infty}_{k=0}\phi_k(r)\Psi_k(\theta),\quad W(r,\theta)=\sum^{\infty}_{k=0}\psi_k(r)\Psi_k(\theta),
    \end{equation}
    where $r=|x|$, $\theta=x/|x|\in \mathbb{S}^{N-1}$, and
    \begin{equation*}
    \phi_k(r)=\int_{\mathbb{S}^{N-1}}v(r,\theta)\Psi_k(\theta)d\theta,\quad \psi_k(r)=\int_{\mathbb{S}^{N-1}}W(r,\theta)\Psi_k(\theta)d\theta.
    \end{equation*}
    Here $\Psi_k(\theta)$ denotes the $k$-th spherical harmonic, i.e., it satisfies
    \begin{equation}\label{deflk}
    -\Delta_{\mathbb{S}^{N-1}}\Psi_k=\lambda_k \Psi_k,
    \end{equation}
    where $\Delta_{\mathbb{S}^{N-1}}$ is the Laplace-Beltrami operator on $\mathbb{S}^{N-1}$ with the standard metric and  $\lambda_k$ is the $k$-th eigenvalue of $-\Delta_{\mathbb{S}^{N-1}}$. We refer to \cite{Mo02} for details.

    It is known that
    \begin{align}\label{Ppwhl2deflklw}
    \Delta (\varphi_k(r)\Psi_k(\theta))
    & = \Psi_k\left(\varphi''_k+\frac{N-1}{r}\varphi'_k\right)
    +\frac{\varphi_k}{r^2}\Delta_{\mathbb{S}^{N-1}}\Psi_k \nonumber\\
    & = \Psi_k\left(\varphi''_k+\frac{N-1}{r}\varphi'_k
    -\frac{\lambda_k}{r^2}\varphi_k\right).
    \end{align}
    Furthermore, it is easy to verify that
    \begin{equation*}
    \frac{\partial (\varphi_k(r)\Psi_k(\theta))}{\partial x_i}=\varphi'_k\frac{x_i}{r}\Psi_k+\varphi_k\frac{\partial\Psi_k}{\partial \theta_l}\frac{\partial\theta_l}{\partial x_i},\quad \mbox{for all}\quad l=1,\ldots,N-1,
    \end{equation*}
    and
    \begin{equation*}
    \sum^{N}_{i=1}\frac{\partial\theta_l}{\partial x_i}x_i=0,\quad \mbox{for all}\quad l=1,\ldots,N-1,
    \end{equation*}
    hence
    \begin{equation}\label{Ppwhl2deflkln}
    \begin{split}
    x\cdot\nabla (\varphi_k(r)\Psi_k(\theta))=\sum^{N}_{i=1}x_i\frac{\partial (\varphi_k(r)\Psi_k(\theta))}{\partial x_i}=\varphi'_kr\Psi_k+\varphi_k\frac{\partial\Psi_k}{\partial \theta_l}\sum^{N}_{i=1}\frac{\partial\theta_l}{\partial x_i}x_i=\varphi'_kr\Psi_k.
    \end{split}
    \end{equation}
    Therefore, by standard regularity theory, putting together (\ref{Ppwhl2deflklw}) and (\ref{Ppwhl2deflkln}) into (\ref{Pwhlp}), the function $(v,W)$ is a solution of (\ref{Pwhlp}) if and only if $(\phi_k,\psi_k)\in \mathcal{C}\times \mathcal{C}$ is a classical solution of the system
    \begin{eqnarray}\label{p2c}
    \left\{ \arraycolsep=1.5pt
       \begin{array}{ll}
        \phi''_k+\frac{N-1+\alpha}{r}\phi'_k
        -\frac{\lambda_k}{r^2}\phi_k+\frac{\psi_k}{r^{\alpha}}=0 \quad \mbox{in}\quad r\in(0,\infty),\\[3mm]
        \psi''_k+\frac{N-1+\alpha}{r}\psi'_k
        -\frac{\lambda_k}{r^2}\psi_k
        +\frac{(2^{**}_{\alpha}-1)C_{N,\alpha}^{2^{**}_{\alpha}-2}}
        {r^{\alpha} (1+r^{2-\alpha})^4}\phi_k=0 \quad \mbox{in}\quad r\in(0,\infty),\\[3mm]
        \phi'_k(0)=\psi'_k(0)=0 \quad\mbox{if}\quad k=0,\quad \mbox{and}\quad \phi_k(0)=\psi_k(0)=0 \quad\mbox{if}\quad k\geq 1,
        \end{array}
    \right.
    \end{eqnarray}
    where
    \[
    \mathcal{C}:=\left\{\omega\in C^2([0,\infty))| \int^\infty_0 \left[\omega''(r)+\frac{N+\alpha-1}{r}\omega'(r)\right]^2
    r^{N+2\alpha-1}\mathrm{d}r<\infty\right\}.
    \]
    Taking the same variation as in the proof of Corollary \ref{thmPbcb}, $|x|=r=s^{1/q}$ where $q=(2-\alpha)/2$ and let
    \begin{equation}\label{p2txy}
    X_k(s)=\phi_k(r),\quad Y_k(s)=q^{-2}\psi_k(r),
    \end{equation}
    that transforms (\ref{p2c}) into the system
    \begin{eqnarray}\label{p2t}
    \left\{ \arraycolsep=1.5pt
       \begin{array}{ll}
        X''_k+\frac{M-1}{s}X'_k-\frac{\lambda_k}{q^2s^2}X_k+Y_k=0 \quad \mbox{in}\quad s\in(0,\infty),\\[3mm]
        Y''_k+\frac{M-1}{s}Y'_k-\frac{\lambda_k}{q^2s^2}Y_k+\frac{(M+4)(M-2)M(M+2)}{(1+s^2)^4}X_k=0 \quad \mbox{in}\quad s\in(0,\infty),\\[3mm]
        X'_k(0)=Y'_k(0)=0 \quad\mbox{if}\quad k=0,\quad \mbox{and}\quad X_k(0)=Y_k(0)=0 \quad\mbox{if}\quad k\geq 1,
        \end{array}
    \right.
    \end{eqnarray}
    in $(X_k,Y_k)\in \widetilde{\mathcal{C}}\times \widetilde{\mathcal{C}}$, where
    \[
    \widetilde{\mathcal{C}}:=\left\{\omega\in C^2([0,\infty))\Big| \int^\infty_0 \left[\omega''(s)+\frac{M-1}{s}\omega'(s)\right]^2 s^{M-1} \mathrm{d}s<\infty\right\},
    \]
    and
    \begin{equation}
    M=\frac{2N}{2-\alpha}>4.
    \end{equation}
    Here we have used the fact
    \begin{small}\begin{equation*}
    q^{-4}(2^{**}_{\alpha}-1)C_{N,\alpha}^{2^{**}_{\alpha}-2}
    =\left[(M-4)(M-2)M(M+2)\right]\left[\frac{2M}{M-4}-1\right]
    =(M+4)(M-2)M(M+2).
    \end{equation*}\end{small}

    Fix $M$ let us now consider the following eigenvalue problem
    \begin{eqnarray}\label{p2te}
    \left\{ \arraycolsep=1.5pt
       \begin{array}{ll}
        X''+\frac{M-1}{s}X'-\frac{\mu}{s^2}X+Y=0 \quad \mbox{in}\quad s\in(0,\infty),\\[3mm]
        Y''+\frac{M-1}{s}Y'-\frac{\mu}{s^2}Y
        +\frac{(M+4)(M-2)M(M+2)}{(1+s^2)^4}X=0 \quad \mbox{in}\quad s\in(0,\infty),
        \end{array}
    \right.
    \end{eqnarray}
    $X\in \widetilde{\mathcal{C}}$. When $M$ is an integer we can study (\ref{p2te}) as the linearized operator of the equation
    \begin{equation*}
    \Delta^2 z=(M-4)(M-2)M(M+2) z^{\frac{M+4}{M-4}},\quad z>0 \quad\mbox{in}\quad \mathbb{R}^M,
    \end{equation*}
    around the standard solution $z(x)=(1+|x|^2)^{-\frac{M-4}{2}}$ (note that we always have $M>4$). In this case, as in \cite[Theorem 2.2]{LW00} (see also \cite{BWW03,Fe02}), we have that
    \begin{equation}\label{ptev}
    \mu_0=0; \quad \mu_1=M-1\quad \mbox{and}\quad X_0(s)=\frac{1-s^2}{(1+s^2)^{\frac{M-2}{2}}}; \quad X_1(s)=\frac{s}{(1+s^2)^{\frac{M-2}{2}}}.
    \end{equation}
    Moreover, even if $M$ is not an integer we readily see that (\ref{ptev}) remains true since \eqref{p2te} is an ODE system and the linearized operator admits at most one positive eigenvalue, see \cite[Proposition 2.4]{BCG21}.
    Here we explain the solutions of \eqref{p2te} with $\mu=0$ must be $X_0$ (up to multiplications). Note that \cite[Lemma 2.4]{BWW03} states that if $X\in \widetilde{\mathcal{C}}$ is a solution of
    \[
    \left[s^{1-M}\frac{\partial}{\partial s}\left(s^{M-1}\frac{\partial}{\partial s}\right)\right]^2X=\nu(1+s^2)^{-4}X \quad \mbox{for}\ \nu>0,
    \]
    with $X(0)=0$, then $X\equiv 0$ (which states that $M$ is an integer, in fact, it also holds for all $M>4$ since the only one step needs to be modified is showing that if $X\not\equiv 0$ then $X$ has only a finite number of positive zeros which requires the conclusions of \cite[Proposition 2]{El77} and \cite[p. 273]{Sw92} and they are indeed correct for all $M>4$). Now, let $\tilde{X}_0\not\equiv 0$ be another solution of \eqref{p2te}. Then $X_0(0), \tilde{X}_0(0)\neq 0$, for otherwise $X_0$ resp. $\tilde{X}_0$ would vanish identically by \cite[Lemma 2.4]{BWW03}. So we can find $\tau\in\mathbb{R}$ such that $X_0(0)=\tau\tilde{X}_0(0)$. But then $X_0-\tau\tilde{X}_0$ also solves \eqref{p2te} and equals zero at origin. By \cite[Lemma 2.4]{BWW03}, one has $X_0-\tau\tilde{X}_0\equiv 0$, and so $\tilde{X}_0$ is a scalar multiple of $X_0$.
    Therefore, we conclude that (\ref{p2t}) has nontrivial solutions if and only if
    \begin{equation*}
    \frac{\lambda_k}{q^2}\in \{0,M-1\},
    \end{equation*}
    where $\lambda_k=k(N-2+k)$, $k\in\mathbb{N}$. If $\lambda_k/q^2=0$ then $k=0$. Moreover, if $\lambda_k/q^2=M-1$,
    that is,
    \begin{equation}\label{ncqf}
    \left(\frac{2-\alpha}{2}\right)^2\left(\frac{2N}{2-\alpha}-1\right)=k(N-2+k),\quad \mbox{for some}\quad k\in\mathbb{N}^+.
    \end{equation}
    Turning back to (\ref{p2c}), if (\ref{ncqf}) holds then
    \begin{equation*}%\label{pye}
    \phi_0(r)=\frac{1-r^{2-\alpha}}{(1+r^{2-\alpha})^{\frac{N-2+\alpha}{2-\alpha}}},\quad
    \phi_k(r)=\frac{r^{\frac{2-\alpha}{2}}}{(1+r^{2-\alpha})^{\frac{N-2+\alpha}{2-\alpha}}},
    \end{equation*}
    for some $k\in\mathbb{N}^+$, otherwise
    \begin{equation*}%\label{pyf}
    \phi_0(r)=\frac{1-r^{2-\alpha}}{(1+r^{2-\alpha})^{\frac{N-2+\alpha}{2-\alpha}}}.
    \end{equation*}
    That is, if (\ref{ncqf}) holds for some $k\in\mathbb{N}^+$, then the space of solutions of (\ref{Pwhlp}) has dimension $1+\frac{(N+2k-2)(N+k-3)!}{(N-2)!k!}$ and is spanned by
    \begin{equation*}
    Z_0(x)=\frac{1-|x|^{2-\alpha}}{(1+|x|^{2-\alpha})^\frac{N-2+\alpha}{2-\alpha}},\quad Z_{k,i}(x)=\frac{|x|^{\frac{2-\alpha}{2}}\Psi_{k,i}(x)}
    {(1+|x|^{2-\alpha})^\frac{N-2+\alpha}{2-\alpha}},
    \end{equation*}
    where $\{\Psi_{k,i}\}$, $i=1,\ldots,\frac{(N+2k-2)(N+k-3)!}{(N-2)!k!}$, form a basis of $\mathbb{Y}_k(\mathbb{R}^N)$, the space of all homogeneous harmonic polynomials of degree $k$ in $\mathbb{R}^N$. Otherwise the space of solutions of (\ref{Pwhlp}) has dimension one and is spanned by $Z_0$, and note that $Z_0\thicksim \frac{\partial U_{\lambda,\alpha}}{\partial \lambda}|_{\lambda=1}$ in this case we say $U$ is non-degenerate. The proof of Theorem \ref{thmpwhl} is now completed.
    \qed

\section{{\bfseries Remainder terms of inequality (\ref{Ppbcm})}}\label{sect:rt}

In this section, we consider the remainder terms of inequality (\ref{Ppbcm}) in radial space $\mathcal{D}^{2,2}_{\alpha,rad}(\mathbb{R}^N)$ and give the proof of Theorem \ref{thmprt}.  We follow the arguments as those in \cite{BE91}, and also \cite{LW00}.

We define $u_\lambda(x):=\lambda^{\frac{N-4+2\alpha}{2}}u(\lambda x)$ for all $\lambda>0$. Thus for simplicity of notations, hereafter of this paper, we write $U_\lambda$ instead of $U_{\lambda,\alpha}$ and $S_\alpha$ instead of $S^{rad}_{N,\alpha}$ if there is no possibility of confusion. Moreover, in order to shorten formulas, for each $u\in \mathcal{D}^{2,2}_{\alpha }(\mathbb{R}^N)$ we denote
    \begin{equation}\label{def:norm}
    \|u\|: =\left(\int_{\mathbb{R}^N} |{\rm div}(|x|^\alpha\nabla u)|^2 \mathrm{d}x\right)^{\frac{1}{2}},
    \quad \|u\|_*: =\left(\int_{\mathbb{R}^N}|u|^{2^{**}_{\alpha}} \mathrm{d}x\right)^{\frac{1}{2^{**}_{\alpha}}}.
    \end{equation}

    Consider the eigenvalue problem
    \begin{equation}\label{Pwhlep}
    \mathrm{div}(|x|^{\alpha}\nabla(\mathrm{div}(|x|^\alpha\nabla v)))=\mu U_{\lambda}^{2^{**}_{\alpha}-2}v \quad \mbox{in}\quad \mathbb{R}^N, \quad v\in \mathcal{D}^{2,2}_\alpha(\mathbb{R}^N).
    \end{equation}
    By a simple scaling argument we know that $\mu$ does not depend on $\lambda$.
    Then from Theorem \ref{thmpwhl} we directly obtain the following conclusion:

    \begin{theorem}\label{propep}
    Suppose $N\geq 2$ and $\frac{4-N}{2}<\alpha<2$. Let $\mu_i$, $i=1,2,\ldots,$ denote the eigenvalues of (\ref{Pwhlep}) in increasing order. Then $\mu_1=1$ is simple with the corresponding eigenfunction $\zeta U_\lambda$ for all $\zeta\in\mathbb{R}\backslash\{0\}$. Moreover,
    \begin{itemize}
    \item[$(i)$]
    if $\alpha$ satisfies (\ref{ncqft}), that is
    \begin{equation*}
    \left(\frac{2-\alpha}{2}\right)^2\left(\frac{2N}{2-\alpha}-1\right)=k(N-2+k),\quad \mbox{for some}\quad k\in\mathbb{N}^+,
    \end{equation*}
    then $\mu_2=\cdots=\mu_{M_k+2}=2^{**}_{\alpha}-1$ with the corresponding $(M_k+1)$-dimensional eigenfunction space spanned by
    \[\left\{\nabla_\lambda U_\lambda,\quad  \lambda^{\frac{N-4+2\alpha}{2}}Z_{k,i}(\lambda x), i=1,\ldots,M_k\right\},\]
    where $M_k=\frac{(N+2k-2)(N+k-3)!}{(N-2)!k!}$, and $Z_{k,i}$ are given as in (\ref{defaezki}).
    \item[$(ii)$]
    otherwise, if (\ref{ncqft}) does not hold, then $\mu_2=2^{**}_{\alpha}-1$ with the corresponding one-dimensional eigenfunction space spanned by $\nabla_\lambda U_\lambda$.
    \end{itemize}
    Furthermore, the definition of eigenvalues indicates $\mu_3>2^{**}_{\alpha}-1$.
    \end{theorem}

    The main ingredient in the proof of Theorem \ref{thmprt} is contained in the lemma below, where the behavior near  $\mathcal{M}_\alpha=\{cU_{\lambda}: c\in\mathbb{R}, \lambda>0\}$ is studied. Here $\mathcal{M}_\alpha$ is the set of extremal functions of inequality \eqref{Ppbcm}.

    \begin{lemma}\label{lemma:rtnm2b}
    Suppose $N\geq 2$ and $\frac{4-N}{2}<\alpha<2$. Then there is $\rho>0$ such that for any sequence $\{u_n\}\subset \mathcal{D}^{2,2}_{\alpha,rad}(\mathbb{R}^N)\backslash \mathcal{M}_\alpha$ satisfying $\inf\limits_{n\in\mathbb{N}}\|u_n\|>0$ and $\mathrm{dist}(u_n,\mathcal{M}_\alpha)\to 0$,
    \begin{equation}\label{rtnmb}
    \liminf\limits_{n\to\infty}\frac{\|u_n\|^2-S_\alpha\|u_n\|^2_*}
    {\mathrm{dist}(u_n,\mathcal{M}_\alpha)^2}\geq \rho.
    \end{equation}
    \end{lemma}

    \begin{proof}
    Let $d_n:={\rm dist}(u_n,\mathcal{M}_\alpha)=\inf\limits_{c\in\mathbb{R}, \lambda>0}\|u_n-cU_\lambda\|\to 0$. When $n$ is sufficiently large, it is easy to verify that there exist $c_n\in\mathbb{R}\backslash \{0\}$ and $\lambda_n>0$ such that $d_n=\|u_n-c_nU_{\lambda_n}\|$, see
    \cite[Lemma 3.2]{DT22}.
    Since $\mathcal{M}_\alpha$ is two-dimensional manifold embedded in $\mathcal{D}^{2,2}_{\alpha,rad}(\mathbb{R}^N)$:
    \[
    (c,\lambda)\in\mathbb{R}\times\mathbb{R}^+\to cU_\lambda\in \mathcal{D}^{2,2}_{\alpha,rad}(\mathbb{R}^N),
    \]
    then the tangential space at $(c_n,\lambda_n)$ is given by
    \[
    T_{c_n U_{\lambda_n}}\mathcal{M}_\alpha={\rm Span}\left\{U_{\lambda_n}, \quad \frac{\partial U_\lambda}{\partial \lambda}\Big|_{\lambda=\lambda_n}\right\}.
    \]
    Anyway we must have that $(u_n-c_n U_{\lambda_n})$ is perpendicular to $T_{c_n U_{\lambda_n}}\mathcal{M}_\alpha$. In particular,
    \begin{align}\label{pul}
    0 & =-\int_{\mathbb{R}^N}|x|^{\alpha}\nabla(\mathrm{div}(|x|^\alpha\nabla U_{\lambda_n}))\cdot\nabla (u_n-c_n U_{\lambda_n})
    \mathrm{d}x \nonumber\\
    & = \int_{\mathbb{R}^N} \mathrm{div}(|x|^\alpha\nabla U_{\lambda_n}) \mathrm{div}(|x|^{\alpha} \nabla (u_n-c_n U_{\lambda_n})) \mathrm{d}x.
    \end{align}
    Furthermore, since we deal with it in radial space, Theorem \ref{propep} implies there exists constant $\mu_3>2^{**}_\alpha-1$ depending only on $N$ and $\alpha$ such that
    \begin{equation}\label{epkeyibbg}
    \mu_3\int_{\mathbb{R}^N}U_{\lambda_n}^{2^{**}_\alpha-2}(u_n-c_n U_{\lambda_n})^2 \mathrm{d}x
    \leq  \int_{\mathbb{R}^N} |{\rm div}(|x|^\alpha\nabla (u_n-c_n U_{\lambda_n}))|^2 \mathrm{d}x.
    \end{equation}
    Let
    \begin{equation}\label{defunwn}
    u_n=c_n U_{\lambda_n}+d_n w_n,
    \end{equation}
     then $w_n$ is perpendicular to $T_{c_n U_{\lambda_n}}\mathcal{M}_\alpha$, $\|w_n\|=1$ and
    \begin{equation*}
    \|u_n\|^2=d_n^2+c_n^2\|U_{\lambda_n}\|^2,
    \end{equation*}
    due to
    \begin{equation*}
    \int_{\mathbb{R}^N} \mathrm{div}(|x|^\alpha\nabla w_n) \mathrm{div}(|x|^{\alpha} \nabla U_{\lambda_n}) \mathrm{d}x
    =\int_{\mathbb{R}^N} \mathrm{div}(|x|^\alpha\nabla U_{\lambda_n}) {\rm div}(|x|^{\alpha} \nabla w_n) \mathrm{d}x=0.
    \end{equation*}
    Then we can rewrite (\ref{epkeyibbg}) as follows:
    \begin{equation}\label{epkeyibbb}
    \int_{\mathbb{R}^N} U_{\lambda_n}^{2^{**}_\alpha-2}w_n^2 \mathrm{d}x\leq \frac{1}{\mu_3}.
    \end{equation}
    By  using Taylor's expansion, it holds that
    \begin{align}\label{epkeyiybb}
    \|u_n\|^{2^{**}_\alpha}_*
    & = \int_{\mathbb{R}^N} |c_n U_{\lambda_n}+d_nw_n|^{2^{**}_\alpha}  \mathrm{d}x \nonumber\\
    & = |c_n|^{2^{**}_\alpha}\int_{\mathbb{R}^N} U_{\lambda_n}^{2^{**}_\alpha} \mathrm{d}x
    +|c_n|^{2^{**}_\alpha-1}d_n 2^{**}_\alpha\int_{\mathbb{R}^N} U_{\lambda_n}^{2^{**}_\alpha-1}w_n \mathrm{d}x  \nonumber\\
    & \quad +\frac{2^{**}_\alpha(2^{**}_\alpha-1)d_n^2  |c_n|^{2^{**}_\alpha-2} }{2}\int_{\mathbb{R}^N} U_{\lambda_n}^{2^{**}_\alpha-2}w_n^2 \mathrm{d}x
    +o(d_n^2)  \nonumber\\
    & = |c_n|^{2^{**}_\alpha}\int_{\mathbb{R}^N} U_{\lambda_n}^{2^{**}_\alpha} \mathrm{d}x
    + \frac{2^{**}_\alpha(2^{**}_\alpha-1)d_n^2  |c_n|^{2^{**}_\alpha-2} }{2} \int_{\mathbb{R}^N} U_{\lambda_n}^{2^{**}_\alpha-2}w_n^2 \mathrm{d}x
    +o(d_n^2),
    \end{align}
    due to
    \begin{equation*}
    \begin{split}
    \int_{\mathbb{R}^N}U_{\lambda_n}^{2^{**}_\alpha-1}w_n \mathrm{d}x
    & = \int_{\mathbb{R}^N}\mathrm{div}(|x|^{\alpha}\nabla({\rm div}(|x|^\alpha\nabla U_{\lambda_n})))w_n \mathrm{d}x \\
    & = \int_{\mathbb{R}^N} \mathrm{div}(|x|^\alpha\nabla U_{\lambda_n}) \mathrm{div}(|x|^{\alpha} \nabla w_n) \mathrm{d}x=0.
    \end{split}
    \end{equation*}
    Then combining with (\ref{epkeyibbb}) and (\ref{epkeyiybb}), by the concavity of $t\mapsto t^{\frac{2}{2^{**}_\alpha}}$ due to $2^{**}_\alpha>2$, we obtain
    \begin{align}\label{epkeyiyxbb}
    \|u_n\|^{2}_*
    & \leq c_n^2\left(\|U_{\lambda_n}\|^{2^{**}_\alpha}_*
    +\frac{2^{**}_\alpha(2^{**}_\alpha-1)d_n^2  c_n^{-2}}{2\mu_3}
    +o(d_n^2)\right)^{\frac{2}{2^{**}_\alpha}}  \nonumber\\
    & = c_n^2\left(\|U_{\lambda_n}\|^2_*+\frac{2}{2^{**}_\alpha}
    \frac{2^{**}_\alpha(2^{**}_\alpha-1)d_n^2  c_n^{-2}}{2\mu_3} \|U_{\lambda_n}\|^{2-2^{**}_\alpha}_*
    +o(d^2)\right) \nonumber\\
    & = c_n^2\|U_{\lambda_n}\|^2_*+ \frac{d_n^2 (2^{**}_\alpha-1)}{\mu_3}\|U_{\lambda_n}\|^{2-2^{**}_\alpha}_*
    +o(d_n^2).
    \end{align}
    Therefore, as $d_n\to 0$ (that is $n$ tends to $\infty$) we have
    \begin{equation*}
    \begin{split}
    \|u_n\|^2-S_\alpha\|u_n\|^{2}_*
    & \geq d_n^2+c_n^2\|U_{\lambda_n}\|^2- S_\alpha\left(c_n^2\|U_{\lambda_n}\|^2_*+ \frac{d_n^2 (2^{**}_\alpha-1)}{\mu_3}\|U_{\lambda_n}\|^{2-2^{**}_\alpha}_*
    +o(d_n^2)\right)  \\
    & = d_n^2 \left(1-\frac{2^{**}_\alpha-1}{\mu_3} S_\alpha \|U_{\lambda_n}\|_*^{2-2^{**}_\alpha}\right)
    +c_n^2(\|U_{\lambda_n}\|^2- S_\alpha\|U_{\lambda_n}\|^2_*)+o(d_n^2)  \\
    & = d_n^2\left(1-\frac{2^{**}_\alpha-1}{\mu_3}\right)+o(d_n^2),
    \end{split}
    \end{equation*}
    since $\|U\|^2=\|U\|_*^{2^{**}_\alpha}
    =S_\alpha^{2^{**}_\alpha/(2^{**}_\alpha-2)}$ and $\|U\|^2=S_\alpha\|U\|_*^2$ for all $U\in\mathcal{M}_\alpha$. Taking $\rho=1-\frac{2^{**}_\alpha-1}{\mu_3}>0$, then (\ref{rtnmb}) follows immediately.
    \end{proof}

\medskip

\noindent {\bf \em Proof of Theorem \ref{thmprt}.}
We argue by contradiction. In fact, if the theorem is false then there exists a sequence $\{u_n\}\subset \mathcal{D}^{2,2}_{\alpha,rad}(\mathbb{R}^N)\backslash \mathcal{M}_\alpha$ such that
    \begin{equation*}
    \frac{\|u_n\|^2-S_\alpha\|u_n\|^2_*}{\mathrm{dist}(u_n,\mathcal{M}_\alpha)^2}\to 0,\quad \mbox{as}\quad n\to \infty.
    \end{equation*}
    By homogeneity, we can assume that $\|u_n\|=1$, and after selecting a subsequence we can assume that $\mathrm{dist}(u_n,\mathcal{M}_\alpha)\to \xi\in[0,1]$ since $\mathrm{dist}(u_n,\mathcal{M}_\alpha)=\inf\limits_{c\in\mathbb{R}, \lambda>0}\|u_n-cU_{\lambda}\|\leq \|u_n\|$. If $\xi=0$, then we have a contradiction by Lemma \ref{lemma:rtnm2b}.

    The other possibility only is that $\xi>0$, that is
    \[\mathrm{dist}(u_n,\mathcal{M}_\alpha)\to \xi>0\quad \mbox{as}\quad n\to \infty,\]
    then we must have
    \begin{equation}\label{wbsi}
    \|u_n\|^2-S_\alpha\|u_n\|^2_*\to 0,\quad \|u_n\|=1.
    \end{equation}
    Since $\{u_n\}\subset \mathcal{D}^{2,2}_{\alpha,rad}(\mathbb{R}^N)\backslash \mathcal{M}_\alpha$ are radial, making the changes that $v_n(s)=u_n(r)$ and $r=s^{2/(2-\alpha)}$, then (\ref{wbsi}) is equivalent to
    \begin{equation}\label{bsiy}
    \int^\infty_0\left[v_n''(s)+\frac{M-1}{s}v_n'(s)\right]^2 s^{M-1}\mathrm{d}s
    -C(M)\left(\int^\infty_0|v_n(s)|^{\frac{2M}{M-4}}s^{M-1}
    \mathrm{d}s\right)^{\frac{M-4}{M}}\to 0,
    \end{equation}
    where $M=\frac{2N}{2-\alpha}>4$ and $C(M)=(M-4)(M-2)M(M+2)\left[\Gamma^2(\frac{M}{2})/(2\Gamma(M))\right]^{\frac{4}{M}}$, see the proof of Corollary \ref{thmPbcb}. When $M$ is an integer,  (\ref{bsiy}) is equivalent to
    \begin{equation}\label{bsib}
    \begin{split}
    \int_{\mathbb{R}^M}|\Delta v_n|^2 \mathrm{d}x
    -S(M)\left(\int_{\mathbb{R}^M}|v_n|^{\frac{2M}{M-4}}
    \mathrm{d}x\right)^{\frac{M-4}{M}}\to 0,
    \end{split}
    \end{equation}
    with
    \[
    \quad \|v_n\|_{\mathcal{D}^{2,2}_{0}(\mathbb{R}^M)}=
    \left(\frac{\omega_{M-1}}{\omega_{N-1}}
    \left(\frac{2}{2-\alpha}\right)^3\right)^{\frac{1}{2}},
    \]
    where $S(M)=\pi^2(M-4)(M-2)M(M+2)\left[\Gamma(\frac{M}{2})/\Gamma(M)\right]^{\frac{4}{M}}$ is the best constant for the embedding of the space $\mathcal{D}^{2,2}_0(\mathbb{R}^M)$ into $L^{2M/(M-4)}(\mathbb{R}^M)$, see \cite{Va93}. In this case, by Lions' concentration and compactness principle (see
    \cite[Theorem \uppercase\expandafter{\romannumeral 1}.1]{Li85-1}), we have that there exists a sequence of positive numbers $\{\tau_n\}$ such that
    \begin{equation*}
    \tau_n^{\frac{M-4}{2}}v_n(\tau_n x)\to V\quad \mbox{in}\quad \mathcal{D}^{2,2}_0(\mathbb{R}^M)\quad \mbox{as}\quad n\to \infty,
    \end{equation*}
    where $V(x)=c(a+|x|^2)^{-(M-4)/2}$ for some $c\neq 0$ and $a>0$, that is
    \begin{equation*}
    \lambda_n^{\frac{N-4+2\alpha}{2}}u_n(\lambda_n x)\to U\quad \mbox{in}\quad \mathcal{D}^{2,2}_\alpha(\mathbb{R}^N)\quad \mbox{as}\quad n\to \infty,
    \end{equation*}
    for some $U\in\mathcal{M}_\alpha$, where $\lambda_n=\tau_n^{\frac{2}{2-\alpha}}$, which implies
    \begin{equation*}
    \mathrm{dist}(u_n,\mathcal{M}_\alpha)=\mathrm{dist}\left(\lambda_n^{\frac{N-4+2\alpha}{2}}u_n(\lambda_n x),\mathcal{M}_\alpha\right)\to 0 \quad \mbox{as}\quad n\to \infty,
    \end{equation*}
    this is a contradiction.

    On the other hand, if $M$ is not an integer, since \eqref{bsiy} is equivalent to
    \begin{equation}\label{bsibni}
    \begin{split}
    \int_0^\infty|s^{1-M}(s^{M-1}v'_n)'|^2 s^{M-1}\mathrm{d}s
    -C(M)\left(\int_0^\infty|v_n|^{\frac{2M}{M-4}}s^{M-1}
    \mathrm{d}s\right)^{\frac{M-4}{M}}\to 0,
    \end{split}
    \end{equation}
    then from the proof of Theorem 1.2 in \cite{dS23} which establishes the Lions' type concentration-compactness principle for high order Sobolev inequality of radial case, we can also get analogous contradiction. Now the proof of Theorem \ref{thmprt} is completed.
    \qed

\section{{\bfseries The finite-dimensional reduction}}\label{sectprp}

    In this section, we consider the perturbation problem (\ref{Pwhp}) and give the proof of Theorem \ref{thmpwhp} by using finite-dimensional reduction method. For simplicity of notations, we write $\int\cdot$ instead of $\int_{\mathbb{R}^N}\cdot \mathrm{d}x$. Hereafter of this section, we always suppose that $N\geq 5$, $0<\alpha<2$ and $h\in L^\infty(\mathbb{R}^N)\cap C(\mathbb{R}^N)$.

    Firstly, we give an equivalent form as the following.

    \begin{proposition}\label{propneq}
    There is $1<C<\infty$ independent of $u$ such that
\begin{align}\label{neq}
\frac{1}{C}\int |\mathrm{div} (|x|^{\alpha}\nabla u)|^2
\leq \int |x|^{2\alpha }|\Delta u|^2
\leq C \int  |\mathrm{div} (|x|^{\alpha}\nabla u)|^2,
\end{align}
for any $u\in C^\infty_0(\mathbb{R}^N)$.
\end{proposition}

\begin{proof}
We follow the arguments as those in \cite[Section 2]{GG22}. Note that
\[
\mathrm{div} (|x|^{\alpha}\nabla u)=|x|^{\alpha}\Delta u+\alpha|x|^{\alpha-2}x\cdot\nabla u,
\]
we see that
\begin{align}\label{neqe}
\int |\mathrm{div} (|x|^{\alpha}\nabla u)|^2
= \int|x|^{2\alpha }|\Delta u|^2
+ 2\alpha \int|x|^{2\alpha-2 }\Delta u(x\cdot\nabla u)
+ \alpha^2 \int|x|^{2\alpha-4 }(x\cdot\nabla u)^2.
\end{align}
The H\"{o}lder inequality implies
\begin{align*}
\left|\int|x|^{2\alpha-2 }\Delta u(x\cdot\nabla u)\right|
& \leq \int|x|^{2\alpha-1 }|\Delta u||\nabla u|
= \int(|x|^{\frac{2\alpha }{2}}|\Delta u|)(|x|^{\frac{2\alpha -2}{2}}|\nabla u|)
\\
& \leq \left(\int|x|^{2\alpha }|\Delta u|^2\right)^{\frac{1}{2}}
\left(\int|x|^{2\alpha -2}|\nabla u|^2\right)^{\frac{1}{2}},
\end{align*}
thus the Young's inequality implies
\begin{align*}
\left|2\alpha \int|x|^{2\alpha-2 }\Delta u(x\cdot\nabla u)\right|
\leq |\alpha|\left(\int|x|^{2\alpha }|\Delta u|^2+\int|x|^{2\alpha -2}|\nabla u|^2\right).
\end{align*}
Furthermore,
\[
\int|x|^{2\alpha-4 }(x\cdot\nabla u)^2
\leq \int|x|^{2\alpha-2 }|\nabla u|^2.
\]
By using the weighted Hardy-Rellich inequality (see \cite{Li86}) we obtain
\[
\int|x|^{2\alpha -2}|\nabla u|^2
\leq D \int|x|^{2\alpha }|\Delta u|^2,
\]
for some $D>0$ independent of $u$, then we deduce the left inequality in \eqref{neq} with $C=1+|\alpha|(1+D)+D\alpha^2$.

Then we show the right inequality in \eqref{neq}. Consider
\begin{align}\label{neqer}
\mathrm{div} (|x|^{\alpha}\nabla u)= v,
\end{align}
we see that
\begin{align}\label{neqer1}
\int |\mathrm{div} (|x|^{\alpha}\nabla u)|^2
=\int|v|^2.
\end{align}
It follows from \eqref{neqer} that
\begin{align*}
|x|^{\alpha }\Delta u=-v-\alpha|x|^{\alpha -2}(x\cdot \nabla u),
\end{align*}
thus
\begin{align}\label{neqer2}
\int|x|^{2\alpha }|\Delta u|^2
=\int|v|^2
+2\alpha\int|x|^{\alpha-2}(x\cdot \nabla u)v
+\alpha^2\int |x|^{2\alpha -4}(x\cdot \nabla u)^2.
\end{align}
Taking $\phi=|x|^{\alpha -2} u$ as a test function to \eqref{neqer2}, we obtain
\[
\int|x|^{\alpha}\nabla u\cdot \nabla \phi
=\int v \phi,
\]
that is,
\[
\int |x|^{2\alpha -2}|\nabla u|^2+(\alpha -2)\int|x|^{2\alpha-\beta-4}
\left(x\cdot \nabla \left(\frac{|u|^2}{2}\right)\right)
=\int |x|^{\alpha -2}vu.
\]
Therefore,
\[
\int |x|^{2\alpha -2}|\nabla u|^2-\frac{\alpha -2}{2}
\int |u|^2\mathrm{div}(|x|^{2\alpha -4}x)
=\int |x|^{\alpha -2}vu,
\]
and
\begin{align}\label{neqere}
\int |x|^{2\alpha-2}|\nabla u|^2+\frac{(2 -\alpha)(N+2\alpha -4)}{2}
\int|x|^{2\alpha -4}|u|^2
=\int |x|^{\alpha -2}vu.
\end{align}
Since $0<\alpha<2$, we obtain that
\[
\frac{(2-\alpha)(N+2\alpha-4)}{2}>0.
\]
It follows from \eqref{neqere} and the H\"{o}lder inequality that
\begin{align*}
\int |x|^{2\alpha-2}|\nabla u|^2
\leq \int |x|^{\alpha-2}vu
\leq \left(\int |v|^2\right)^{\frac{1}{2}}\left(\int |x|^{2\alpha-4}|u|^2\right)^{\frac{1}{2}}.
\end{align*}
By using the weighted Hardy inequality (see \cite{CKN84}) we obtain
\[
\int|x|^{2\alpha-4}|u|^2
\leq E \int|x|^{2\alpha-2}|\nabla u|^2,
\]
for some $E>0$ independent of $u$, we see that
\begin{align*}
\int |x|^{2\alpha-2}|\nabla u|^2
\leq E\int |v|^2.
\end{align*}
It follows from \eqref{neqer1} and \eqref{neqer2} that
\begin{align*}
\int|x|^{2\alpha}|\Delta u|^2
& \leq \int|v|^2
+|\alpha|\left(\int|x|^{2\alpha-2}|\nabla u|^2
+\int|v|^2
\right)
+\alpha^2\int |x|^{2\alpha-2}|\nabla u|^2
\\
& = (1+|\alpha|)\int|v|^2 +(|\alpha|+\alpha^2)\int |x|^{2\alpha-2}|\nabla u|^2
\\
& \leq \left[1+|\alpha|+E(|\alpha|+\alpha^2)\right]
\int |\mathrm{div} (|x|^{\alpha}\nabla u)|^2.
\end{align*}
Therefore, we obtain the right inequality in \eqref{neq} with $C=1+|\alpha|+E(|\alpha|+\alpha^2)$.

Note that the constant $C$ in \eqref{neq} can be chosen to
\[
C=\max\{1+|\alpha|(1+D)+D\alpha^2, 1+|\alpha|(1+E)+E\alpha^2\}.
\]
Now the proof is completed.
\end{proof}

In \cite[Lemma A.1]{MS14}, Musina and Sreenadh established the following inequality (we only mention the second-order case): let $N\geq 5$ and $2\leq q\leq \frac{2N}{N-4}$, if
\[
\mu_2(a):=\inf_{u\in C^2_0(\mathbb{R}^N)\setminus\{0\},\ u(x)=u(|x|)}\frac{\int |x|^{a}|\Delta u|^2}{\int |x|^{a-4}|u|^2}>0,
\]
then there exists a constant $C>0$ such that
\begin{equation}\label{cknccg}
    \int |x|^{a}|\Delta u|^2  \geq C \left(\int |x|^{-N+q\frac{N-4+a}{2}}|u|^{q}\right)^{\frac{2}{q}} ,\quad
    \mbox{for all}\quad u\in C^2_0(\mathbb{R}^N).
    \end{equation}
    Note that when $N\geq 5$ and $0<\alpha<2$, we have $2<2^{**}_\alpha<\frac{2N}{N-4}$ and $\mu_2(2\alpha)=\frac{(N-4+2\alpha)(N-2\alpha)}{4}>0$ (see \cite{MS12}), then from Proposition \ref{propneq} and \eqref{cknccg} we deduce that
    \begin{equation*}
    \int |\mathrm{div}(|x|^\alpha\nabla u)|^2
    \geq C_1\int |x|^{2\alpha}|\Delta u|^2
    \geq C_2 \left(\int |u|^{2^{**}_\alpha}\right)^{\frac{2}{2^{**}_\alpha}} ,\quad
    \mbox{for all}\quad u\in C^\infty_0(\mathbb{R}^N),
    \end{equation*}
Since $\mathcal{D}^{2,2}_{\alpha}(\mathbb{R}^N)$ is dense in $C^\infty_0(\mathbb{R}^N)$,  we deduce that $\mathcal{D}^{2,2}_{\alpha}(\mathbb{R}^N)\hookrightarrow
    L^{2^{**}_\alpha}(\mathbb{R}^N)$ continuously.

    Now, let us begin to deal with the perturbation  problem \eqref{Pwhp}. Given $h\in L^\infty(\mathbb{R}^N)\cap C(\mathbb{R}^N)$, we put
    \begin{equation}\label{defH}
    H[u]=\frac{1}{2^{**}_\alpha}\int h(x) u^{2^{**}_\alpha}_+.
    \end{equation}
    For $\varepsilon\in\mathbb{R}$, we introduce the perturbed energy functional $\mathcal{J}_\varepsilon$ and also the unperturbed energy functional $\mathcal{J}_0$ given by
    \begin{align}\label{defje}
    \mathcal{J}_\varepsilon[u]
    =   \mathcal{J}_0[u]-\varepsilon H[u]
    =  \frac{1}{2}\int
    |\mathrm{div}(|x|^{\alpha}\nabla u)|^2
    -\frac{1}{2^{**}_\alpha}\int (1+\varepsilon h(x))u^{2^{**}_\alpha}_+,
    \end{align}
    for $u\in \mathcal{D}^{2,2}_{\alpha}(\mathbb{R}^N)$.
    Evidently, $\mathcal{J}_\varepsilon\in C^2$ and any critical point $u$ of $\mathcal{J}_\varepsilon$ is a weak solution to
    \begin{equation*}
    \mathrm{div}(|x|^{\alpha}\nabla(\mathrm{div}(|x|^\alpha\nabla u)))=(1+\varepsilon h(x))|x|^{-\alpha} u^{2^{**}_\alpha-1}_+.
    \end{equation*}
    If $u\neq 0$ and $|\varepsilon|\|h\|_\infty< 1$, then $u$ is positive by the strong maximum principle. Hence, $u$ solves (\ref{Pwhp}).

    Define now
    \begin{equation}\label{deful}
    \mathcal{U}:=\left\{U_{\lambda}(x)
    =\lambda^{\frac{N-4+2\alpha}{2}}U_{1}(\lambda x):  \lambda>0\right\},
    \end{equation}
    where $U_{\lambda}:=U_{\lambda,\alpha}$ as in \eqref{defula}, and for each $\lambda>0$ we define
    \begin{align}\label{defeip0}
    \mathcal{E}_\lambda:=\left\{\omega\in \mathcal{D}^{2,2}_{\alpha}(\mathbb{R}^N): \left\langle\omega,\frac{\partial U_{\lambda}}{\partial \lambda}\right\rangle_\alpha
    =0\right\}.
    \end{align}
    Here $\left\langle\omega,\frac{\partial U_{\lambda}}{\partial \lambda}\right\rangle_\alpha
    =\int\mathrm{div}(|x|^\alpha\nabla \omega)\mathrm{div}
    \left(|x|^\alpha\nabla \frac{\partial U_{\lambda}}{\partial \lambda}\right)$.
    From Theorem \ref{thmpwhl}, under the assumption $N\geq 5$ and $0<\alpha<2$ we know that the manifold $\mathcal{U}$ is non-degenerate. Then taking the same argument as in \cite[Corollary 3.2]{FS03}, it is easy to verify that $\mathcal{J}''_0[U_1]$ is a self-adjoint Fredholm operator of index zero which maps the space $\mathcal{D}^{2,2}_\alpha(\mathbb{R}^N)$ into $T_{U_1}\mathcal{U}^\perp$, and $\mathcal{J}''_0[U_1]\in \mathfrak{L}(T_{U_1}\mathcal{U}^\perp)$ is invertible.

    \begin{lemma}\label{lemhw}
    There exists a constant $C_1=C_1(N,\alpha,\|h\|_\infty)>0$ such that for any $\lambda>0$, $\omega\in \mathcal{D}^{2,2}_{\alpha}(\mathbb{R}^N)$, it holds that
    \begin{equation}\label{gh0}
    |H[U_{\lambda}+\omega]|\leq C_1(\||h|^{\frac{1}{2^{**}_\alpha}}U_{\lambda}\|^{2^{**}_\alpha}_*
    +\|\omega\|^{2^{**}_\alpha}),
    \end{equation}
    \begin{equation}\label{gh1}
    \|H'[U_{\lambda}+\omega]\|\leq C_1(\||h|^{\frac{1}{2^{**}_\alpha}}U_{\lambda}\|^{2^{**}_\alpha-1}_*
    +\|\omega\|^{2^{**}_\alpha-1}),
    \end{equation}
    \begin{equation}\label{gh2}
    \|H''[U_{\lambda}+\omega]\|\leq C_1(\||h|^{\frac{1}{2^{**}_\alpha}}U_{\lambda}\|^{2^{**}_\alpha-2}_*
    +\|\omega\|^{2^{**}_\alpha-2}).
    \end{equation}
    Moreover, if $\lim\limits_{|x|\to 0}h(x)=\lim\limits_{|x|\to \infty}h(x)=0$ then
    \begin{equation}\label{ghu}
    \||h|^{\frac{1}{2^{**}_\alpha}}U_{\lambda}\|_*\to 0 \quad\mbox{as}\quad \lambda\to 0\quad\mbox{or}\quad \lambda\to \infty.
    \end{equation}
    See the definitions of $\|\cdot\|$ and $\|\cdot\|_*$ as in (\ref{def:norm}).
    \end{lemma}

    \begin{proof}
    We will only show (\ref{gh2}) as (\ref{gh0})-(\ref{gh1}) follow analogously. By H\"{o}lder's inequality and since the embedding $\mathcal{D}^{2,2}_{\alpha}(\mathbb{R}^N)\hookrightarrow L^{2^{**}_\alpha}(\mathbb{R}^N)$ is continuous, we have
    \begin{equation*}
    \begin{split}
    \|H''[U_{\lambda}+\omega]\|
    & \leq (2^{**}_\alpha-1)\sup_{\|g_1\|,\|g_1\|\leq 1} \int |h(x)||U_{\lambda}+\omega|^{2^{**}_\alpha-2}
    |g_1||g_2| \\
    & \leq (2^{**}_\alpha-1)\|h\|^{\frac{2}{2^{**}_\alpha}}_\infty \sup_{\|g_1\|,\|g_1\|\leq 1} \||h|^{\frac{1}{2^{**}_\alpha}}
    (U_{\lambda}+\omega)\|^{2^{**}_\alpha-2}_*
    \|g_1\|_*\|g_2\|_* \\
    & \leq c(N,\alpha,\|h\|_\infty)\||h|^{\frac{1}{2^{**}_\alpha}}
    (U_{\lambda}+\omega)\|^{2^{**}_\alpha-2}_*.
    \end{split}
    \end{equation*}
    Then by using the triangle inequality and again $\mathcal{D}^{2,2}_{\alpha}(\mathbb{R}^N)\hookrightarrow L^{2^{**}_\alpha}(\mathbb{R}^N)$ continuously, we can directly obtain (\ref{gh2}).

    Under the additional assumption $\lim\limits_{|x|\to 0}h(x)=\lim\limits_{|x|\to \infty}h(x)=0$, (\ref{ghu}) follows by the dominated convergence theorem and
    \begin{equation*}
    \int |h(x)|U_{\lambda}^{2^{**}_\alpha}
    =\int |h(\lambda^{-1} x)|U_1^{2^{**}_\alpha} .
    \end{equation*}
    \end{proof}

    In order to deal with the problem $\mathcal{J}'_\varepsilon[u]=0$ for $\varepsilon$ close to zero, we combine variational methods with the standard Lyapunov-Schmit reduction.

    \begin{lemma}\label{lemreg}
    There exists a small constant $\varepsilon_0>0$ such that for each $\lambda>0$ and $|\varepsilon|<\varepsilon_0$, there is a smooth function $\omega(\lambda,\varepsilon):
    (0,\infty)\times(-\varepsilon_0,\varepsilon_0)
    \to \mathcal{D}^{2,2}_{\alpha}(\mathbb{R}^N)
    $ satisfying
    \begin{equation}\label{we}
    \omega(\lambda,\varepsilon)\in \mathcal{E}_\lambda,
    \end{equation}
    \begin{equation}\label{jes}
    \mathcal{J}'_\varepsilon[U_\lambda
    +\omega(\lambda,\varepsilon)]\eta=0,
    \quad \mbox{for all}\ \eta\in \mathcal{E}_\lambda,
    \end{equation}
    \begin{equation}\label{wgx}
    \|\omega(\lambda,\varepsilon)\|\leq C_2|\varepsilon|,
    \end{equation}
    for some constant $C>0$.
    Moreover, if $\lim\limits_{|x|\to 0}h(x)=\lim\limits_{|x|\to \infty}h(x)=0$ then
    \begin{equation}\label{wgt0}
    \|\omega(\lambda,\varepsilon)\|\to 0 \quad\mbox{as}\quad \lambda\to 0\quad\mbox{or}\quad \lambda\to \infty,
    \end{equation}
    uniformly with respect to $\varepsilon$.
    \end{lemma}

    \begin{proof}
    The proof can be deduced directly from the proof of \cite[Lemma 4.2]{DT22} with minor changes, so we omit it.
    \end{proof}

    Under the assumptions of Lemma \ref{lemreg}, for each $\varepsilon$ satisfying $|\varepsilon|<\varepsilon_0$ we define
    \begin{equation}\label{defule}
    \mathcal{U^\varepsilon}:=\left\{u\in \mathcal{D}^{2,2}_{\alpha}(\mathbb{R}^N)|
    u=U_{\lambda}+\omega(\lambda,\varepsilon),\quad \lambda\in(0,\infty)\right\},
    \end{equation}
    where $\varepsilon_0>0$ and $\omega(\lambda,\varepsilon)\in \mathcal{E}_\lambda$ are obtained in Lemma \ref{lemreg}. Note that $\mathcal{U^\varepsilon}$ is a one-dimensional manifold. The next lemma will show that finding critical points for functional can be reduced to a finite dimensional problem.

    \begin{lemma}\label{lemcuve}
    Under the assumptions of Lemma \ref{lemreg}, then there is $\varepsilon_1\in (0,\varepsilon_0)$ such that for each $|\varepsilon|<\varepsilon_1$ the manifold $\mathcal{U^\varepsilon}$ is a natural constraint for $\mathcal{J}_\varepsilon$, that is, every critical point of $\mathcal{J}_\varepsilon|_{\mathcal{U^\varepsilon}}$ is the critical point of $\mathcal{J}_\varepsilon$.
    \end{lemma}

    \begin{proof}
    The proof can be deduced directly from the proof of \cite[Lemma 4.3]{DT22} with minor changes, so we omit it.
    \end{proof}

    Now, we are in position to prove the main result.

\medskip

    \noindent {\bf \em Proof of Theorem \ref{thmpwhp}.}
    Choosing $|\varepsilon|$ sufficiently small, then Lemma \ref{lemreg} indicates for each $\lambda>0$ there is $\omega(\lambda,\varepsilon)$ satisfying \eqref{we}-\eqref{wgt0}.

    Let $u^\varepsilon_\lambda=U_\lambda+\omega(\lambda,\varepsilon)$ then
    \begin{equation*}
    \mathcal{J}_\varepsilon[u^\varepsilon_\lambda]
    =\mathcal{J}_0[U_{\lambda}]
    +\frac{1}{2}(\|u^\varepsilon_\lambda\|^2-\|U_{\lambda}\|^2)
    -\frac{1}{2^{**}_\alpha}\int
    (1+\varepsilon h(x))((u^\varepsilon_\lambda)^{2^{**}_\alpha}_+ -U_{\lambda}^{2^{**}_\alpha})
    -\varepsilon H[U_{\lambda}].
    \end{equation*}
    Recall that $\|U_{\lambda}\|=\|U_{1}\|$ does not depend on $\lambda$, and (\ref{wgx}) indicates $\|u^\varepsilon_\lambda-U_{\lambda}\|
    =\|\omega(\lambda,\varepsilon)\|=o(1)$ as $\lambda\to 0$ or $\lambda\to \infty$ uniformly with respect to $\varepsilon$, therefore
    \begin{equation*}
    |\|u^\varepsilon_\lambda\|^2-\|U_{\lambda}\|^2|
    \leq (\|\omega(\lambda,\varepsilon)\|+2\|U_{\lambda}\|)\|u^\varepsilon_\lambda-U_{\lambda}\|=o(1).
    \end{equation*}
    Moreover, by H\"{o}lder's inequality and $\mathcal{D}^{2,2}_{\alpha}(\mathbb{R}^N)\hookrightarrow L^{2^{**}_\alpha}(\mathbb{R}^N)$ continuously we obtain
    \begin{equation*}
    \begin{split}
    \left|\int (1+\varepsilon h(x))((u^\varepsilon_\lambda)^{2^{**}_\alpha}_+ -U_{\lambda}^{2^{**}_\alpha}) \right|
    & \leq C\int ((u^\varepsilon_\lambda)^{2^{**}_\alpha-1}_+ +U_{\lambda}^{2^{**}_\alpha-1})|(u^\varepsilon_\lambda)_+ -U_{\lambda}|
    \\
    & \leq  C\|u^\varepsilon_\lambda-U_{\lambda}\|=o(1).
    \end{split}
    \end{equation*}
    Finally, combining with (\ref{ghu}) which indicates $|H[U_{\lambda}]|=o(1)$ as $\lambda\to 0$ or $\lambda\to \infty$, we conclude that
    \begin{equation*}
    \Gamma_{\varepsilon}(\lambda):=\mathcal{J}_\varepsilon[u^\varepsilon_\lambda]=\mathcal{J}_0[U_{\lambda}]+o(1)=\mathcal{J}_0[U_{1}]+o(1),\quad \mbox{as}\quad \lambda\to 0\quad \mbox{or} \quad \lambda\to\infty,
    \end{equation*}
    that is,
    \begin{equation*}
    \lim_{\lambda\to 0}\Gamma_{\varepsilon}(\lambda)=\lim_{\lambda\to \infty}\Gamma_{\varepsilon}(\lambda)=\mathcal{J}_0[U_1],
    \end{equation*}
    uniformly with respect to $\varepsilon$.
    Thus, $\Gamma_{\varepsilon}$ has at least one critical point $\lambda_{\varepsilon}$. Therefore by Lemma \ref{lemcuve},  $u_\varepsilon:=u^\varepsilon_{\lambda_\varepsilon}$ is a critical point for $\mathcal{J}_\varepsilon$, furthermore \eqref{wgx} indicates $\|\omega(\lambda_\varepsilon,\varepsilon)\|=O(\varepsilon)$. Now, the proof of Theorem \ref{thmpwhp} is completed.
    \qed

\medskip

\noindent{\bfseries Acknowledgements}

The authors wish to thank the referee for his/her very helpful comments and suggestions which indeed improve the quality of the paper.

The research has been supported by National Natural Science Foundation of China (No. 12371121 and 11971392), and Chongqing Graduate Student Research Innovation Project (No. CYB23107).

    \end{document}